# An Enumeration of Certain Equivalence Classes of Directed Multigraphs Having 0 to 5 Nodes with Each Node Having 2 Incoming and 2 Outgoing Arcs

C. C. Briggs
*Center for Academic Computing, Penn State University, University Park, PA 16802*
March 22, 1999

**Abstract.** An enumeration is given of certain equivalence classes of directed multigraphs having 0 to 5 nodes with each node having 2 incoming and 2 outgoing arcs.

This paper presents an enumeration of certain equivalence classes of directed multigraphs having 0 to 5 nodes with each node having 2 incoming and 2 outgoing arcs.

The classes—which consist of distinct but isomorphic multigraphs—are (1) designated by sets of 3 numbers, the $1^{st}$, $2^{nd}$, and $3^{rd}$ of which specify the numbers of nodes of the component multigraphs, the ranks, and the cardinalities of the classes, respectively, and (2) typified by representative multigraphs characterized by pictures and by corresponding polynomials in variables individually signifying the occurrences of directed arcs, each variable having 2 subscripts, the $1^{st}$ and $2^{nd}$ of which specify the nodes from which and at which an arc emanates and terminates, respectively. (For 0 nodes, the representative multigraph is the null multigraph, the picture for which is intentionally left blank and the polynomial taken to be equal to 1.)

For a check, note that, for $p$ nodes, the cardinalities of the classes add up to $(2p)!/2^p$.

Some numerical features of the classes appear in Table 1 (see below).

TABLE 1. SOME NUMERICAL FEATURES OF THE CLASSES

| NUMBER OF NODES | NUMBER OF CLASSES | SUM OF CARDINALITIES |
|---|---|---|
| 0 | **1** | 1 |
| 1 | **1** | 1 |
| 2 | **3** | 6 |
| 3 | **8** | 90 |
| 4 | **25** | 2520 |
| 5 | **85** | 113,400 |

0 NODES

0,1
Class: 0,1,1
Polynomial: 1

1 NODE

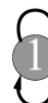

1,1,1

1,1
Class: 1,1,1
Polynomial: $x_{11} x_{11}$

2 NODES

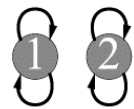

2,1,1

2,1
Class: 2,1,1
Polynomial: $x_{11} x_{11} x_{22} x_{22}$

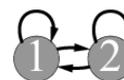

2,2,4

2,2
Class: 2,2,4
Polynomial: $x_{11} x_{12} x_{22} x_{21}$

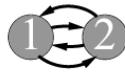

2,3
Class: 2,3,1
Polynomial: $x_{12}\, x_{12}\, x_{21}\, x_{21}$

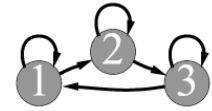

3,4
Class: 3,4,16
Polynomial: $x_{11}\, x_{12}\, x_{22}\, x_{23}\, x_{33}\, x_{31}$

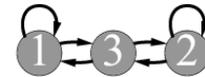

3,5
Class: 3,5,24
Polynomial: $x_{11}\, x_{13}\, x_{22}\, x_{23}\, x_{31}\, x_{32}$

# 3 Nodes

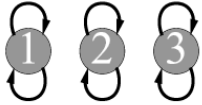

3,1
Class: 3,1,1
Polynomial: $x_{11}\, x_{11}\, x_{22}\, x_{22}\, x_{33}\, x_{33}$

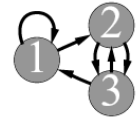

3,6
Class: 3,6,24
Polynomial: $x_{11}\, x_{12}\, x_{23}\, x_{23}\, x_{31}\, x_{32}$

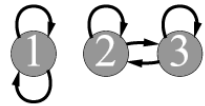

3,2
Class: 3,2,12
Polynomial: $x_{11}\, x_{11}\, x_{22}\, x_{23}\, x_{33}\, x_{32}$

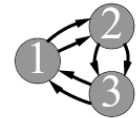

3,7
Class: 3,7,2
Polynomial: $x_{12}\, x_{12}\, x_{23}\, x_{23}\, x_{31}\, x_{31}$

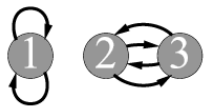

3,3
Class: 3,3,3
Polynomial: $x_{11}\, x_{11}\, x_{23}\, x_{23}\, x_{32}\, x_{32}$

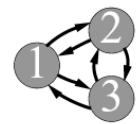

3,8
Class: 3,8,8
Polynomial: $x_{12}\, x_{13}\, x_{21}\, x_{23}\, x_{31}\, x_{32}$

# 4 NODES

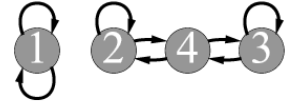

4,5
Class: 4,5,96
Polynomial: $x_{11}\ x_{11}\ x_{22}\ x_{24}\ x_{33}\ x_{34}\ x_{42}\ x_{43}$

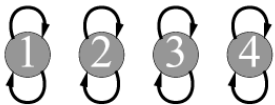

4,1
Class: 4,1,1
Polynomial: $x_{11}\ x_{11}\ x_{22}\ x_{22}\ x_{33}\ x_{33}\ x_{44}\ x_{44}$

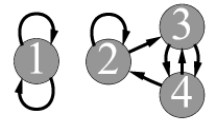

4,6
Class: 4,6,96
Polynomial: $x_{11}\ x_{11}\ x_{22}\ x_{23}\ x_{34}\ x_{34}\ x_{42}\ x_{43}$

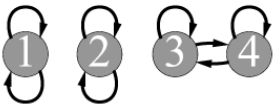

4,2
Class: 4,2,24
Polynomial: $x_{11}\ x_{11}\ x_{22}\ x_{22}\ x_{33}\ x_{34}\ x_{44}\ x_{43}$

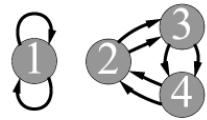

4,7
Class: 4,7,8
Polynomial: $x_{11}\ x_{11}\ x_{23}\ x_{23}\ x_{34}\ x_{34}\ x_{42}\ x_{42}$

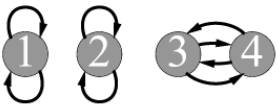

4,3
Class: 4,3,6
Polynomial: $x_{11}\ x_{11}\ x_{22}\ x_{22}\ x_{34}\ x_{34}\ x_{43}\ x_{43}$

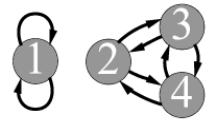

4,8
Class: 4,8,32
Polynomial: $x_{11}\ x_{11}\ x_{23}\ x_{24}\ x_{32}\ x_{34}\ x_{42}\ x_{43}$

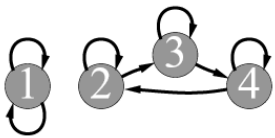

4,4
Class: 4,4,64
Polynomial: $x_{11}\ x_{11}\ x_{22}\ x_{23}\ x_{33}\ x_{34}\ x_{44}\ x_{42}$

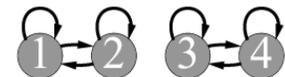

4,9
Class: 4,9,48
Polynomial: $x_{11}\ x_{12}\ x_{22}\ x_{21}\ x_{33}\ x_{34}\ x_{44}\ x_{43}$

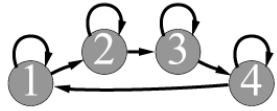

4,10
Class: 4,10,96
Polynomial: $x_{11}\, x_{12}\, x_{22}\, x_{23}\, x_{33}\, x_{34}\, x_{44}\, x_{41}$

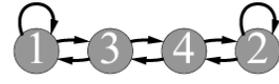

4,15
Class: 4,15,192
Polynomial: $x_{11}\, x_{13}\, x_{22}\, x_{24}\, x_{31}\, x_{34}\, x_{42}\, x_{43}$

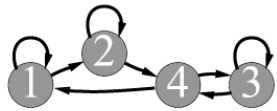

4,11
Class: 4,11,384
Polynomial: $x_{11}\, x_{12}\, x_{22}\, x_{24}\, x_{33}\, x_{34}\, x_{41}\, x_{43}$

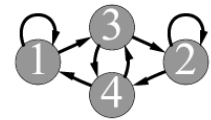

4,16
Class: 4,16,192
Polynomial: $x_{11}\, x_{13}\, x_{22}\, x_{24}\, x_{32}\, x_{34}\, x_{41}\, x_{43}$

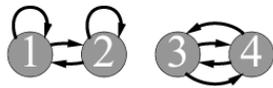

4,12
Class: 4,12,24
Polynomial: $x_{11}\, x_{12}\, x_{22}\, x_{21}\, x_{34}\, x_{34}\, x_{43}\, x_{43}$

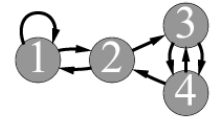

4,17
Class: 4,17,192
Polynomial: $x_{11}\, x_{12}\, x_{21}\, x_{23}\, x_{34}\, x_{34}\, x_{42}\, x_{43}$

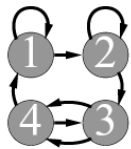

4,13
Class: 4,13,192
Polynomial: $x_{11}\, x_{12}\, x_{22}\, x_{23}\, x_{34}\, x_{34}\, x_{41}\, x_{43}$

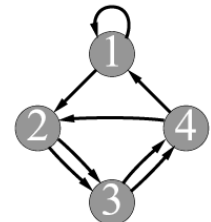

4,18
Class: 4,18,96
Polynomial: $x_{11}\, x_{12}\, x_{23}\, x_{23}\, x_{34}\, x_{34}\, x_{41}\, x_{42}$

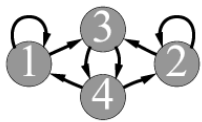

4,14
Class: 4,14,96
Polynomial: $x_{11}\, x_{13}\, x_{22}\, x_{23}\, x_{34}\, x_{34}\, x_{41}\, x_{42}$

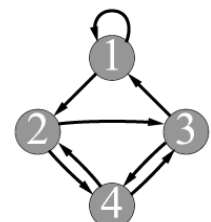

4,19
Class: 4,19,384
Polynomial: $x_{11}\, x_{12}\, x_{23}\, x_{24}\, x_{31}\, x_{34}\, x_{42}\, x_{43}$

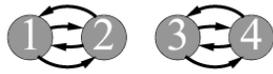

4,20,3

**4,20**
Class: 4,20,3
Polynomial: $x_{12}\ x_{12}\ x_{21}\ x_{21}\ x_{34}\ x_{34}\ x_{43}\ x_{43}$

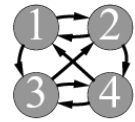

4,25,96

**4,25**
Class: 4,25,96
Polynomial: $x_{12}\ x_{13}\ x_{21}\ x_{24}\ x_{32}\ x_{34}\ x_{41}\ x_{43}$

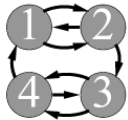

4,21,48

**4,21**
Class: 4,21,48
Polynomial: $x_{12}\ x_{12}\ x_{21}\ x_{23}\ x_{34}\ x_{34}\ x_{41}\ x_{43}$

## 5 NODES

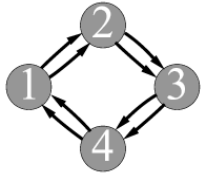

4,22,6

**4,22**
Class: 4,22,6
Polynomial: $x_{12}\ x_{12}\ x_{23}\ x_{23}\ x_{34}\ x_{34}\ x_{41}\ x_{41}$

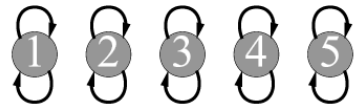

5,1,1

**5,1**
Class: 5,1,1
Polynomial: $x_{11}\ x_{11}\ x_{22}\ x_{22}\ x_{33}\ x_{33}\ x_{44}\ x_{44}\ x_{55}\ x_{55}$

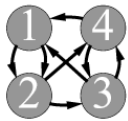

4,23,96

**4,23**
Class: 4,23,96
Polynomial: $x_{12}\ x_{12}\ x_{23}\ x_{24}\ x_{31}\ x_{34}\ x_{41}\ x_{43}$

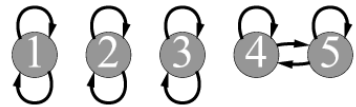

5,2,40

**5,2**
Class: 5,2,40
Polynomial: $x_{11}\ x_{11}\ x_{22}\ x_{22}\ x_{33}\ x_{33}\ x_{44}\ x_{45}\ x_{55}\ x_{54}$

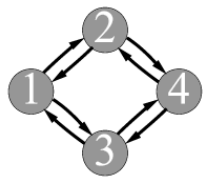

4,24,48

**4,24**
Class: 4,24,48
Polynomial: $x_{12}\ x_{13}\ x_{21}\ x_{24}\ x_{31}\ x_{34}\ x_{42}\ x_{43}$

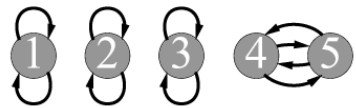

5,3,10

**5,3**
Class: 5,3,10
Polynomial: $x_{11}\ x_{11}\ x_{22}\ x_{22}\ x_{33}\ x_{33}\ x_{45}\ x_{45}\ x_{54}\ x_{54}$

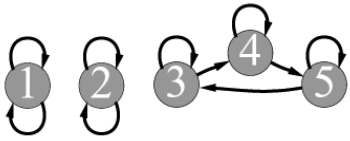

5,4
Class: 5,4,160
Polynomial: $x_{11}\,x_{11}\,x_{22}\,x_{22}\,x_{33}\,x_{34}\,x_{44}\,x_{45}\,x_{55}\,x_{53}$

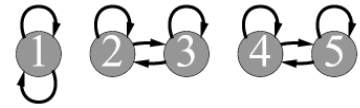

5,9
Class: 5,9,240
Polynomial: $x_{11}\,x_{11}\,x_{22}\,x_{23}\,x_{33}\,x_{32}\,x_{44}\,x_{45}\,x_{55}\,x_{54}$

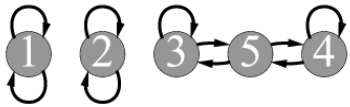

5,5
Class: 5,5,240
Polynomial: $x_{11}\,x_{11}\,x_{22}\,x_{22}\,x_{33}\,x_{35}\,x_{44}\,x_{45}\,x_{53}\,x_{54}$

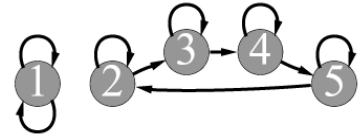

5,10
Class: 5,10,480
Polynomial: $x_{11}\,x_{11}\,x_{22}\,x_{23}\,x_{33}\,x_{34}\,x_{44}\,x_{45}\,x_{55}\,x_{52}$

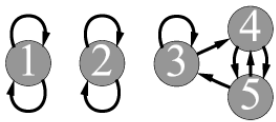

5,6
Class: 5,6,240
Polynomial: $x_{11}\,x_{11}\,x_{22}\,x_{22}\,x_{33}\,x_{34}\,x_{45}\,x_{45}\,x_{53}\,x_{54}$

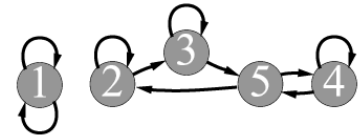

5,11
Class: 5,11,1920
Polynomial: $x_{11}\,x_{11}\,x_{22}\,x_{23}\,x_{33}\,x_{35}\,x_{44}\,x_{45}\,x_{52}\,x_{54}$

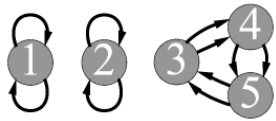

5,7
Class: 5,7,20
Polynomial: $x_{11}\,x_{11}\,x_{22}\,x_{22}\,x_{34}\,x_{34}\,x_{45}\,x_{45}\,x_{53}\,x_{53}$

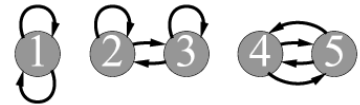

5,12
Class: 5,12,120
Polynomial: $x_{11}\,x_{11}\,x_{22}\,x_{23}\,x_{33}\,x_{32}\,x_{45}\,x_{45}\,x_{54}\,x_{54}$

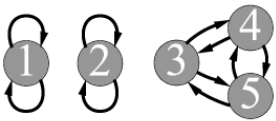

5,8
Class: 5,8,80
Polynomial: $x_{11}\,x_{11}\,x_{22}\,x_{22}\,x_{34}\,x_{35}\,x_{43}\,x_{45}\,x_{53}\,x_{54}$

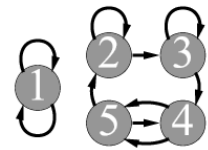

5,13
Class: 5,13,960
Polynomial: $x_{11}\,x_{11}\,x_{22}\,x_{23}\,x_{33}\,x_{34}\,x_{45}\,x_{45}\,x_{52}\,x_{54}$

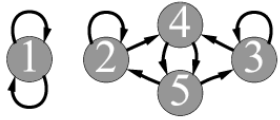

5,14
Class: 5,14,480
Polynomial: $x_{11}\ x_{11}\ x_{22}\ x_{24}\ x_{33}\ x_{34}\ x_{45}\ x_{45}\ x_{52}\ x_{53}$

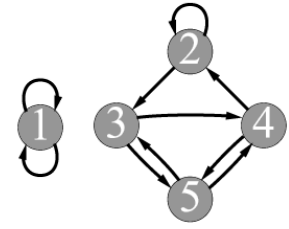

5,19
Class: 5,19,1920
Polynomial: $x_{11}\ x_{11}\ x_{22}\ x_{23}\ x_{34}\ x_{35}\ x_{42}\ x_{45}\ x_{53}\ x_{54}$

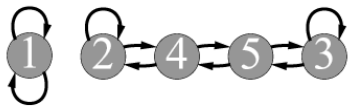

5,15
Class: 5,15,960
Polynomial: $x_{11}\ x_{11}\ x_{22}\ x_{24}\ x_{33}\ x_{35}\ x_{42}\ x_{45}\ x_{53}\ x_{54}$

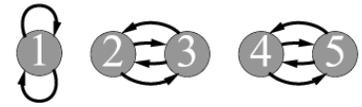

5,20
Class: 5,20,15
Polynomial: $x_{11}\ x_{11}\ x_{23}\ x_{23}\ x_{32}\ x_{32}\ x_{45}\ x_{45}\ x_{54}\ x_{54}$

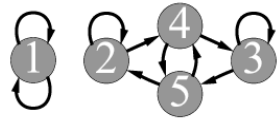

5,16
Class: 5,16,960
Polynomial: $x_{11}\ x_{11}\ x_{22}\ x_{24}\ x_{33}\ x_{35}\ x_{43}\ x_{45}\ x_{52}\ x_{54}$

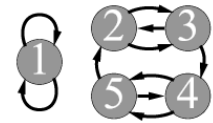

5,21
Class: 5,21,240
Polynomial: $x_{11}\ x_{11}\ x_{23}\ x_{23}\ x_{32}\ x_{34}\ x_{45}\ x_{45}\ x_{52}\ x_{54}$

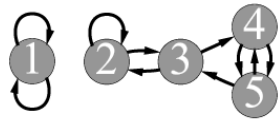

5,17
Class: 5,17,960
Polynomial: $x_{11}\ x_{11}\ x_{22}\ x_{23}\ x_{32}\ x_{34}\ x_{45}\ x_{45}\ x_{53}\ x_{54}$

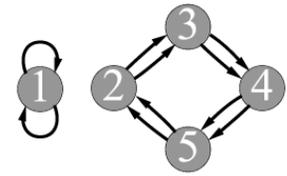

5,22
Class: 5,22,30
Polynomial: $x_{11}\ x_{11}\ x_{23}\ x_{23}\ x_{34}\ x_{34}\ x_{45}\ x_{45}\ x_{52}\ x_{52}$

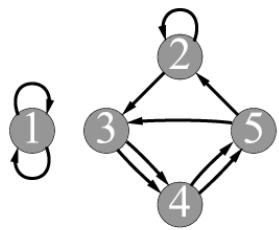

5,18
Class: 5,18,480
Polynomial: $x_{11}\ x_{11}\ x_{22}\ x_{23}\ x_{34}\ x_{34}\ x_{45}\ x_{45}\ x_{52}\ x_{53}$

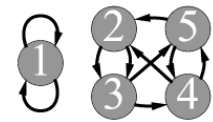

5,23
Class: 5,23,480
Polynomial: $x_{11}\ x_{11}\ x_{23}\ x_{23}\ x_{34}\ x_{35}\ x_{42}\ x_{45}\ x_{52}\ x_{54}$

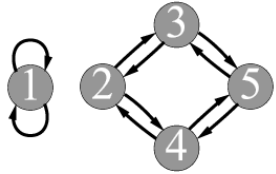

5,24,240

5,24
Class: 5,24,240
Polynomial: $x_{11}\, x_{11}\, x_{23}\, x_{24}\, x_{32}\, x_{35}\, x_{42}\, x_{45}\, x_{53}\, x_{54}$

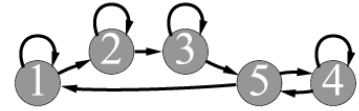

5,29,3840

5,29
Class: 5,29,3840
Polynomial: $x_{11}\, x_{12}\, x_{22}\, x_{23}\, x_{33}\, x_{35}\, x_{44}\, x_{45}\, x_{51}\, x_{54}$

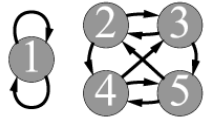

5,25,480

5,25
Class: 5,25,480
Polynomial: $x_{11}\, x_{11}\, x_{23}\, x_{24}\, x_{32}\, x_{35}\, x_{43}\, x_{45}\, x_{52}\, x_{54}$

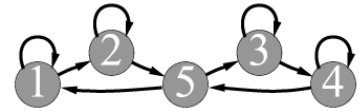

5,30,1920

5,30
Class: 5,30,1920
Polynomial: $x_{11}\, x_{12}\, x_{22}\, x_{25}\, x_{33}\, x_{34}\, x_{44}\, x_{45}\, x_{51}\, x_{53}$

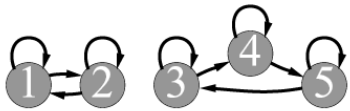

5,26,640

5,26
Class: 5,26,640
Polynomial: $x_{11}\, x_{12}\, x_{22}\, x_{21}\, x_{33}\, x_{34}\, x_{44}\, x_{45}\, x_{55}\, x_{53}$

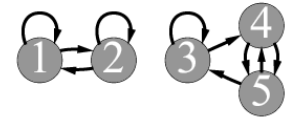

5,31,960

5,31
Class: 5,31,960
Polynomial: $x_{11}\, x_{12}\, x_{22}\, x_{21}\, x_{33}\, x_{34}\, x_{45}\, x_{45}\, x_{53}\, x_{54}$

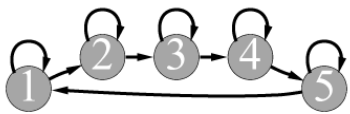

5,27,768

5,27
Class: 5,27,768
Polynomial: $x_{11}\, x_{12}\, x_{22}\, x_{23}\, x_{33}\, x_{34}\, x_{44}\, x_{45}\, x_{55}\, x_{51}$

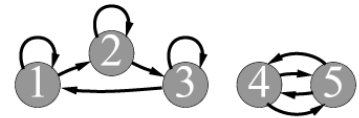

5,32,160

5,32
Class: 5,32,160
Polynomial: $x_{11}\, x_{12}\, x_{22}\, x_{23}\, x_{33}\, x_{31}\, x_{45}\, x_{45}\, x_{54}\, x_{54}$

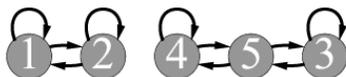

5,28,960

5,28
Class: 5,28,960
Polynomial: $x_{11}\, x_{12}\, x_{22}\, x_{21}\, x_{33}\, x_{35}\, x_{44}\, x_{45}\, x_{53}\, x_{54}$

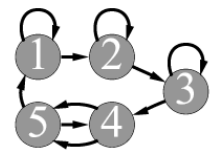

5,33,1920

5,33
Class: 5,33,1920
Polynomial: $x_{11}\, x_{12}\, x_{22}\, x_{23}\, x_{33}\, x_{34}\, x_{45}\, x_{45}\, x_{51}\, x_{54}$

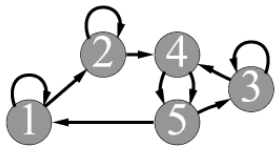

5,34,1920

5,34
Class: 5,34,1920
Polynomial: $x_{11}\,x_{12}\,x_{22}\,x_{24}\,x_{33}\,x_{34}\,x_{45}\,x_{45}\,x_{51}\,x_{53}$

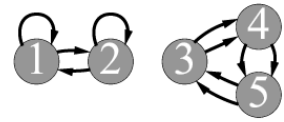

5,39,80

5,39
Class: 5,39,80
Polynomial: $x_{11}\,x_{12}\,x_{22}\,x_{21}\,x_{34}\,x_{34}\,x_{45}\,x_{45}\,x_{53}\,x_{53}$

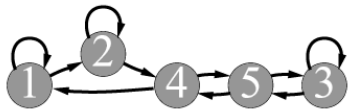

5,35,3840

5,35
Class: 5,35,3840
Polynomial: $x_{11}\,x_{12}\,x_{22}\,x_{24}\,x_{33}\,x_{35}\,x_{41}\,x_{45}\,x_{53}\,x_{54}$

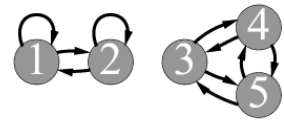

5,40,320

5,40
Class: 5,40,320
Polynomial: $x_{11}\,x_{12}\,x_{22}\,x_{21}\,x_{34}\,x_{35}\,x_{43}\,x_{45}\,x_{53}\,x_{54}$

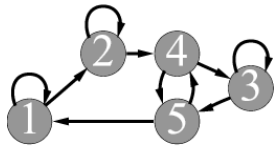

5,36,3840

5,36
Class: 5,36,3840
Polynomial: $x_{11}\,x_{12}\,x_{22}\,x_{24}\,x_{33}\,x_{35}\,x_{43}\,x_{45}\,x_{51}\,x_{54}$

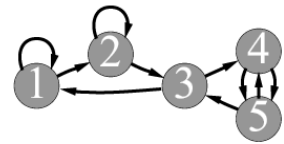

5,41,1920

5,41
Class: 5,41,1920
Polynomial: $x_{11}\,x_{12}\,x_{22}\,x_{23}\,x_{31}\,x_{34}\,x_{45}\,x_{45}\,x_{53}\,x_{54}$

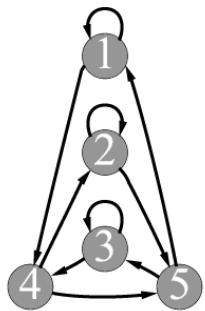

5,37,1920

5,37
Class: 5,37,1920
Polynomial: $x_{11}\,x_{14}\,x_{22}\,x_{25}\,x_{33}\,x_{34}\,x_{42}\,x_{45}\,x_{51}\,x_{53}$

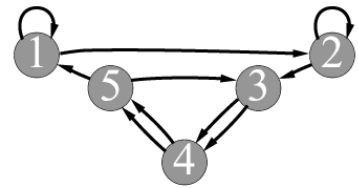

5,42,960

5,42
Class: 5,42,960
Polynomial: $x_{11}\,x_{12}\,x_{22}\,x_{23}\,x_{34}\,x_{34}\,x_{45}\,x_{45}\,x_{51}\,x_{53}$

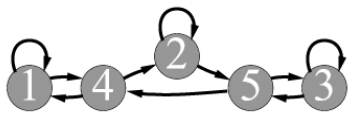

5,38,3840

5,38
Class: 5,38,3840
Polynomial: $x_{11}\,x_{14}\,x_{22}\,x_{25}\,x_{33}\,x_{35}\,x_{41}\,x_{42}\,x_{53}\,x_{54}$

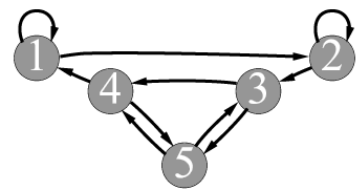

5,43,3840

5,43
Class: 5,43,3840
Polynomial: $x_{11}\,x_{12}\,x_{22}\,x_{23}\,x_{34}\,x_{35}\,x_{41}\,x_{45}\,x_{53}\,x_{54}$

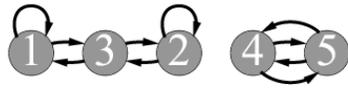

5,44,240

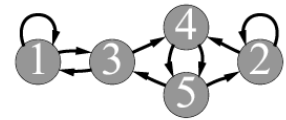

5,49,1920

**5,44**
Class: 5,44,240
Polynomial: $x_{11}\, x_{13}\, x_{22}\, x_{23}\, x_{31}\, x_{32}\, x_{45}\, x_{45}\, x_{54}\, x_{54}$

**5,49**
Class: 5,49,1920
Polynomial: $x_{11}\, x_{13}\, x_{22}\, x_{24}\, x_{31}\, x_{34}\, x_{45}\, x_{45}\, x_{52}\, x_{53}$

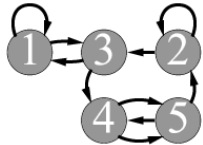

5,45,1920

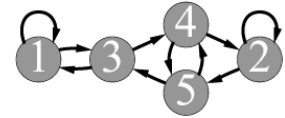

5,50,3840

**5,45**
Class: 5,45,1920
Polynomial: $x_{11}\, x_{13}\, x_{22}\, x_{23}\, x_{31}\, x_{34}\, x_{45}\, x_{45}\, x_{52}\, x_{54}$

**5,50**
Class: 5,50,3840
Polynomial: $x_{11}\, x_{13}\, x_{22}\, x_{25}\, x_{31}\, x_{34}\, x_{42}\, x_{45}\, x_{53}\, x_{54}$

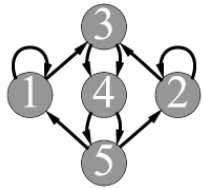

5,46,480

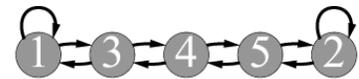

5,51,1920

**5,46**
Class: 5,46,480
Polynomial: $x_{11}\, x_{13}\, x_{22}\, x_{23}\, x_{34}\, x_{34}\, x_{45}\, x_{45}\, x_{51}\, x_{52}$

**5,51**
Class: 5,51,1920
Polynomial: $x_{11}\, x_{13}\, x_{22}\, x_{25}\, x_{31}\, x_{34}\, x_{43}\, x_{45}\, x_{52}\, x_{54}$

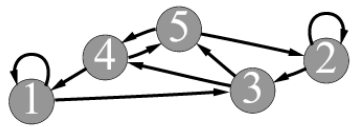

5,47,1920

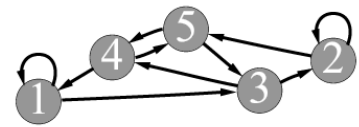

5,52,3840

**5,47**
Class: 5,47,1920
Polynomial: $x_{11}\, x_{13}\, x_{22}\, x_{23}\, x_{34}\, x_{35}\, x_{41}\, x_{45}\, x_{52}\, x_{54}$

**5,52**
Class: 5,52,3840
Polynomial: $x_{11}\, x_{13}\, x_{22}\, x_{25}\, x_{32}\, x_{34}\, x_{41}\, x_{45}\, x_{53}\, x_{54}$

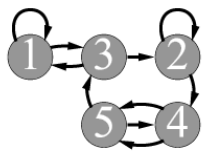

5,48,1920

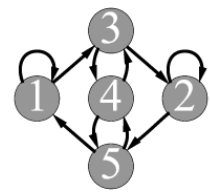

5,53,1920

**5,48**
Class: 5,48,1920
Polynomial: $x_{11}\, x_{13}\, x_{22}\, x_{24}\, x_{31}\, x_{32}\, x_{45}\, x_{45}\, x_{53}\, x_{54}$

**5,53**
Class: 5,53,1920
Polynomial: $x_{11}\, x_{13}\, x_{22}\, x_{25}\, x_{32}\, x_{34}\, x_{43}\, x_{45}\, x_{51}\, x_{54}$

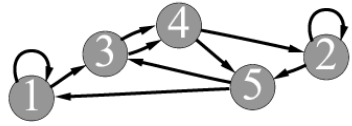

5,54    5,54,1920
Class: 5,54,1920
Polynomial: $x_{11}\,x_{13}\,x_{22}\,x_{25}\,x_{34}\,x_{34}\,x_{42}\,x_{45}\,x_{51}\,x_{53}$

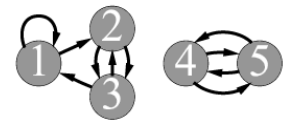

5,59    5,59,240
Class: 5,59,240
Polynomial: $x_{11}\,x_{12}\,x_{23}\,x_{23}\,x_{31}\,x_{32}\,x_{45}\,x_{45}\,x_{54}\,x_{54}$

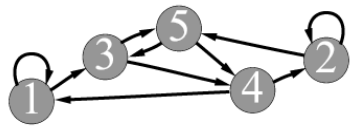

5,55    5,55,1920
Class: 5,55,1920
Polynomial: $x_{11}\,x_{13}\,x_{22}\,x_{25}\,x_{34}\,x_{35}\,x_{41}\,x_{42}\,x_{53}\,x_{54}$

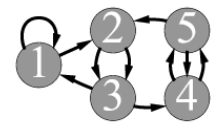

5,60    5,60,960
Class: 5,60,960
Polynomial: $x_{11}\,x_{12}\,x_{23}\,x_{23}\,x_{31}\,x_{34}\,x_{45}\,x_{45}\,x_{52}\,x_{54}$

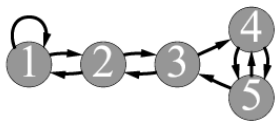

5,56    5,56,1920
Class: 5,56,1920
Polynomial: $x_{11}\,x_{12}\,x_{21}\,x_{23}\,x_{32}\,x_{34}\,x_{45}\,x_{45}\,x_{53}\,x_{54}$

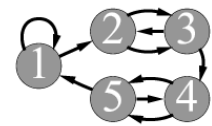

5,61    5,61,960
Class: 5,61,960
Polynomial: $x_{11}\,x_{12}\,x_{23}\,x_{23}\,x_{32}\,x_{34}\,x_{45}\,x_{45}\,x_{51}\,x_{54}$

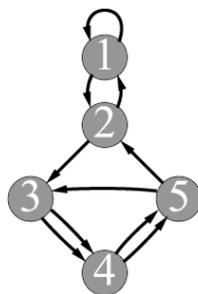

5,57    5,57,960
Class: 5,57,960
Polynomial: $x_{11}\,x_{12}\,x_{21}\,x_{23}\,x_{34}\,x_{34}\,x_{45}\,x_{45}\,x_{52}\,x_{53}$

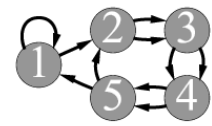

5,62    5,62,480
Class: 5,62,480
Polynomial: $x_{11}\,x_{12}\,x_{23}\,x_{23}\,x_{34}\,x_{34}\,x_{45}\,x_{45}\,x_{51}\,x_{52}$

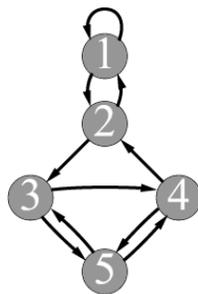

5,58    5,58,3840
Class: 5,58,3840
Polynomial: $x_{11}\,x_{12}\,x_{21}\,x_{23}\,x_{34}\,x_{35}\,x_{42}\,x_{45}\,x_{53}\,x_{54}$

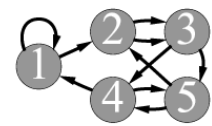

5,63    5,63,1920
Class: 5,63,1920
Polynomial: $x_{11}\,x_{12}\,x_{23}\,x_{23}\,x_{34}\,x_{35}\,x_{41}\,x_{45}\,x_{52}\,x_{54}$

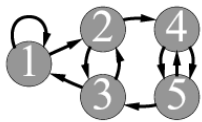

5,64,1920

**5,64**
Class: 5,64,1920
Polynomial: $x_{11}\, x_{12}\, x_{23}\, x_{24}\, x_{31}\, x_{32}\, x_{45}\, x_{45}\, x_{53}\, x_{54}$

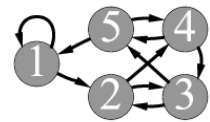

5,69,3840

**5,69**
Class: 5,69,3840
Polynomial: $x_{11}\, x_{12}\, x_{23}\, x_{24}\, x_{32}\, x_{35}\, x_{43}\, x_{45}\, x_{51}\, x_{54}$

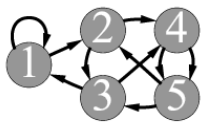

5,65,1920

**5,65**
Class: 5,65,1920
Polynomial: $x_{11}\, x_{12}\, x_{23}\, x_{24}\, x_{31}\, x_{34}\, x_{45}\, x_{45}\, x_{52}\, x_{53}$

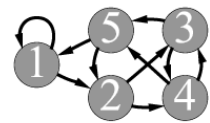

5,70,1920

**5,70**
Class: 5,70,1920
Polynomial: $x_{11}\, x_{12}\, x_{23}\, x_{24}\, x_{34}\, x_{35}\, x_{43}\, x_{45}\, x_{51}\, x_{52}$

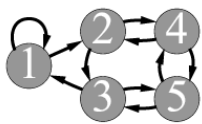

5,66,3840

**5,66**
Class: 5,66,3840
Polynomial: $x_{11}\, x_{12}\, x_{23}\, x_{24}\, x_{31}\, x_{35}\, x_{42}\, x_{45}\, x_{53}\, x_{54}$

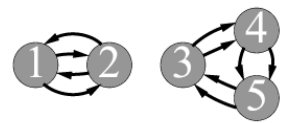

5,71,20

**5,71**
Class: 5,71,20
Polynomial: $x_{12}\, x_{12}\, x_{21}\, x_{21}\, x_{34}\, x_{34}\, x_{45}\, x_{45}\, x_{53}\, x_{53}$

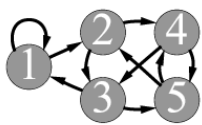

5,67,3840

**5,67**
Class: 5,67,3840
Polynomial: $x_{11}\, x_{12}\, x_{23}\, x_{24}\, x_{31}\, x_{35}\, x_{43}\, x_{45}\, x_{52}\, x_{54}$

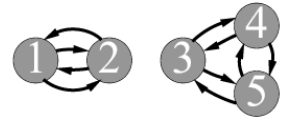

5,72,80

**5,72**
Class: 5,72,80
Polynomial: $x_{12}\, x_{12}\, x_{21}\, x_{21}\, x_{34}\, x_{35}\, x_{43}\, x_{45}\, x_{53}\, x_{54}$

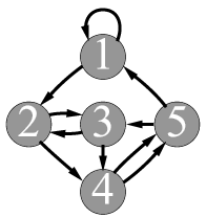

5,68,1920

**5,68**
Class: 5,68,1920
Polynomial: $x_{11}\, x_{12}\, x_{23}\, x_{24}\, x_{32}\, x_{34}\, x_{45}\, x_{45}\, x_{51}\, x_{53}$

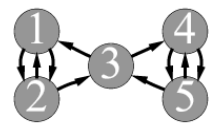

5,73,480

**5,73**
Class: 5,73,480
Polynomial: $x_{12}\, x_{12}\, x_{21}\, x_{23}\, x_{31}\, x_{34}\, x_{45}\, x_{45}\, x_{53}\, x_{54}$

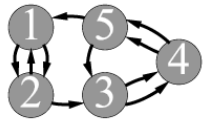

5,74
Class: 5,74,480
Polynomial: $x_{12}\,x_{12}\,x_{21}\,x_{23}\,x_{34}\,x_{34}\,x_{45}\,x_{45}\,x_{51}\,x_{53}$

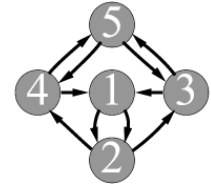

5,79
Class: 5,79,960
Polynomial: $x_{12}\,x_{12}\,x_{23}\,x_{24}\,x_{31}\,x_{35}\,x_{41}\,x_{45}\,x_{53}\,x_{54}$

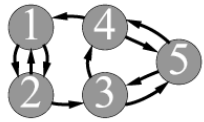

5,75
Class: 5,75,1920
Polynomial: $x_{12}\,x_{12}\,x_{21}\,x_{23}\,x_{34}\,x_{35}\,x_{41}\,x_{45}\,x_{53}\,x_{54}$

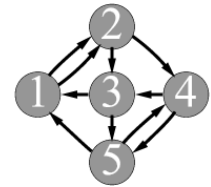

5,80
Class: 5,80,1920
Polynomial: $x_{12}\,x_{12}\,x_{23}\,x_{24}\,x_{31}\,x_{35}\,x_{43}\,x_{45}\,x_{51}\,x_{54}$

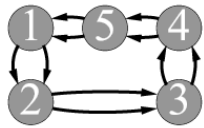

5,76
Class: 5,76,24
Polynomial: $x_{12}\,x_{12}\,x_{23}\,x_{23}\,x_{34}\,x_{34}\,x_{45}\,x_{45}\,x_{51}\,x_{51}$

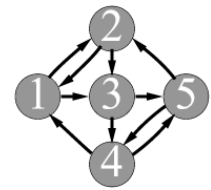

5,81
Class: 5,81,1920
Polynomial: $x_{12}\,x_{13}\,x_{21}\,x_{23}\,x_{34}\,x_{35}\,x_{41}\,x_{45}\,x_{52}\,x_{54}$

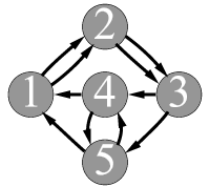

5,77
Class: 5,77,480
Polynomial: $x_{12}\,x_{12}\,x_{23}\,x_{23}\,x_{34}\,x_{35}\,x_{41}\,x_{45}\,x_{51}\,x_{54}$

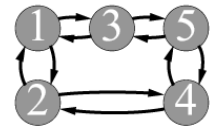

5,82
Class: 5,82,384
Polynomial: $x_{12}\,x_{13}\,x_{21}\,x_{24}\,x_{31}\,x_{35}\,x_{42}\,x_{45}\,x_{53}\,x_{54}$

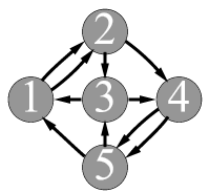

5,78
Class: 5,78,480
Polynomial: $x_{12}\,x_{12}\,x_{23}\,x_{24}\,x_{31}\,x_{34}\,x_{45}\,x_{45}\,x_{51}\,x_{53}$

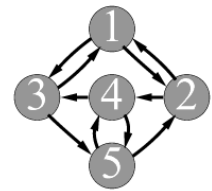

5,83
Class: 5,83,1920
Polynomial: $x_{12}\,x_{13}\,x_{21}\,x_{24}\,x_{31}\,x_{35}\,x_{43}\,x_{45}\,x_{52}\,x_{54}$

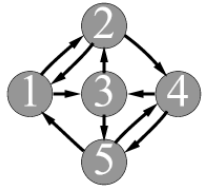

5,84
Class: 5,84,1920
Polynomial: $x_{12}\, x_{13}\, x_{21}\, x_{24}\, x_{32}\, x_{35}\, x_{43}\, x_{45}\, x_{51}\, x_{54}$

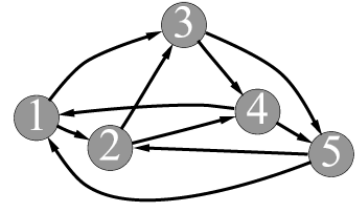

5,85
Class: 5,85,768
Polynomial: $x_{12}\, x_{13}\, x_{23}\, x_{24}\, x_{34}\, x_{35}\, x_{41}\, x_{45}\, x_{51}\, x_{52}$